\documentclass[11pt]{article}

\input{amssym.def}
\input{amssym}
\usepackage{eufrak}

\newtheorem{theorem}{Theorem}

\newtheorem{lemma}[theorem]{Lemma}
\newtheorem{proposition}[theorem]{Proposition}
\newtheorem{remark}[theorem]{Remark}

\def\hl{\hat{l}}
\def\hL{\hat{L}}

\def\un{\underline}
\def\qq{q^{-1}}

\def\Tr{\mathrm{Tr}}
\def\Trr{\Tr_R}

\def\ZZ{{\cal{Z}}}
\def\OO{{\cal{O}}_\mu}
\def\vm{V_{\mm}}

\def\hLL{\hat{{\cal L}}}
\def\LL{{\cal{L}}}
\def\MM{{\cal{M}}}

\def\de{\delta}

\def\ot{\otimes}
\def\C{{\Bbb C}}
\def\Sym{{\rm Sym\, }}
\def\vv{V^{\otimes 2}}
\def\RR{R^{-1}}

\def\mm{{\bf m}}
\def\ov{\overline}
\def\hmu{\hat{\mu}}
\def\hnu{\hat{\hnu}}

\def\la{{\lambda}}
\def\al{{\alpha}}

\def\be{\begin{equation}}
\def\ee{\end{equation}}

\begin{document}

\makeatletter
\renewcommand{\theequation}{{\thesection}.{\arabic{equation}}}
\@addtoreset{equation}{section} \makeatother

\title{Braided Gelfand-Zetlin algebras and their semiclassical counterparts}

\author{\rule{0pt}{7mm} Dimitry
Gurevich\thanks{gurevich@ihes.fr}\\
{\small\it Higher School of Modern Mathematics,}\\
{\small\it Moscow Institute of Physics and Technology}\\
\rule{0pt}{7mm} Pavel Saponov\thanks{Pavel.Saponov@ihep.ru}\\
{\small\it
National Research University Higher School of Economics,}\\
{\small\it 20 Myasnitskaya Ulitsa, Moscow 101000, Russian Federation}\\
{\small\it and}\\
{\small\it Institute for High Energy Physics, NRC "Kurchatov Institute"}\\
{\small \it Protvino 142281, Russian Federation}}

\maketitle

\begin{abstract}
We construct analogs of the Gelfand-Zetlin algebras in  the Reflection Equation  algebras, corresponding to Hecke symmetries, mainly  to those coming from the
Quantum Groups $U_q(sl(N))$. Corresponding semiclassical (i.e. Poisson) counterparts of the Gelfand-Zetlin algebras are described.
\end{abstract}

{\bf AMS Mathematics Subject Classification, 2020:} 17B37, 81R50


\section{Introduction}

Irreducible finite-dimensional  $U(gl(N))$-modules are labeled by  patterns
$$
\mm=(\mm_1\geq \mm_2 \geq\dots \geq \mm_N),
$$
where $\mm_i$ are non-negative integers. The module, corresponding to the pattern $\mm$, will be denoted by $\vm$.
 It is known that the center $Z(U(gl(N)))$ of the algebra $U(gl(N))$ is generated by the elements $p_k=\Tr L^k$, $1\le k\le N$,
where the {\em generating matrix} $L=\|l_i^j\|_{1\leq i,j\leq N}$ is composed of the entries $l_i^j$ subject to the commutation relations
$$
[l_i^j, l_k^r]=l_i^r\de_k^j-l_k^j\de_i^r.
$$
The elements $p_k$ are called {\em power sums} or {\em Gelfand elements}.

Let  $U(gl(N-1))$ be a subalgebra of $U(gl(N))$, generated by the elements $l_i^j$ with indexes $1\leq i,j \leq N-1$. Then a $U(gl(N))$-module $\vm$
can be decomposed into the direct sum of irreducible $U(gl(N-1))$-modules. In the seminal paper  \cite{GZ} the authors  claimed: a $U(gl(N-1))$-module $V_{\mm'}$
is included in the decomposition of the module $\vm$ iff its pattern $\mm'=(\mm_1'\dots \mm_{N-1}')$ subject to the following {\em interlacing condition}:
$$
\mm_1\geq \mm_1'\geq   \mm_2\geq \mm_2' \geq\dots \mm_{N-1}'    \geq \mm_N.
$$
and its multiplicity is 1.

By continuing this procedure, one gets a chain of algebras
\be
U(gl(N))\supset U(gl(N-1))\supset\dots \supset U(gl(2))\supset U(gl(1))
\label{chain}
\ee
and a set of total patterns of the form:
\be
\MM=\left(
\begin{array}{ccccccccc}
\mm_1&&\mm_2&&\mm_3&&\dots&&\mm_N\\
&\mm_1'&&\mm_2'&&\dots&&\mm_{N-1}' &\\
&&\mm_1''&&\dots&&\mm_{N-2}''&&\\
&&&&\dots&&&&\\
\end{array}\right)
\label{pat}
\ee
where any two neighbouring rows are subject  to the interlacing condition.

The lowest element of each total pattern $\MM$ enumerates a one-dimensional subspace of the vector space $\vm$. The directional vector $e_\MM$ of this subspace is a common
eigenvector of elements of all centers $Z(U(gl(k))$,  $1\leq k \leq N$. Moreover, the set of vectors $e_\MM$ spans the whole space $\vm$.

The union of all centers $Z(U(gl(k))$ (denoted $\ZZ$) is a commutative subalgebra of $U(gl(N))$. The vectors $e_\MM$ are common eigenvectors for all elements of $\ZZ$.
The subalgebra $\ZZ$ and the basis composed of the elements $e_\MM$ are called Gelfand-Zetlin (GZ) ones.

Passing from the $gl(N)$ Lie bracket to the corresponding Poisson one $\{\,,\,\}_{gl(N)}$, defined on the commutative algebra   $\Sym(gl(N))$, we get  a semiclassical  counterpart
of the GZ algebra. Its elements are classical analogs of the elements from  $\ZZ$. This classical version of the GZ algebras  was studies in  \cite{GuS, KW}. Since the Poisson bracket
$\{\,,\,\}_{gl(N)}$ admits a restriction to any $GL(N)$-orbit in $\Sym(gl(N))^*$ (where it is called the Kirillov-Kostant-Souriau bracket), it is also possible to restrict the GZ algebra to such an orbit. We are mainly interested in the restrictions to generic orbits\footnote{We call an orbit  {\em semisimple} if it is  the orbit of a diagonal matrix. We call a semisimple orbit
{\em generic}, if the eigenvalues of the corresponding matrix are pairwise distinct.}.

There are known numerous generalizations of the GZ construction for other algebras and modules. Thus, analogs of the GZ basis are constructed in modules over other simple
Lie algebras, some super-algebras and  infinite-dimensional algebras.  In addition, in some publications  there are considered analogs of GZ modules over quantum groups 
and  Yangians. We refer the reader to the paper \cite{M} for numerous references. We only mention the pioneering paper \cite{UTS}, where
a GZ type basis was constructed in $U_q(gl(N))$-modules.

The main objective of the current paper is to construct braided analogs of the  GZ algebras. Let us precise that the term {\em braided} stands for  Reflection Equation (RE) algebras
and all related objects. Any RE algebra $\LL(R)$, we are dealing with, is associated with a Hecke symmetry $R$ which completely defines the 
structure of $\LL(R)$. The best known Hecke symmetries are the Drinfeld-Jimbo ones
coming from the quantum groups $U_q(sl(N))$. They are given by formula  (\ref{He}) below. These Hecke symmetries and all corresponding objects will be called {\em standard} or {\it super-standard} if a Hecke symmetry $R$ comes from a super quantum group $U_q(sl(m|n))$.
However, the family of Hecke symmetries and, consequently, the family of  RE algebras is much larger.

 If an RE algebra can be embedded into another one similarly to
the embedding 
$$U(gl(N-1))\hookrightarrow  U(gl(N))$$
 described above, the union of their centers is a commutative algebra. We call it a {\em  braided GZ type algebra}. For some
Hecke symmetries it is possible to construct a braided version of the chain (\ref{chain}) (an analog of full flags).  Such are standard, super-standard symmetries and some other close to them.

If a Hecke symmetry $R$ defining the  largest RE algebra of the chain is a deformation of the usual flip $P$, then there exists the
corresponding classical $r$-matrix which is subject to the classical Yang-Baxter equation. The related RE algebra $\LL(R)$ has a
semiclassical counterpart, i.e. a commutative algebra $\Sym(gl(N))$ endowed with a Poisson bracket, denoted $\{\,,\,\}_r$. This Poisson bracket is compatible with the linear Poisson bracket $\{\,,\,\}_{gl(N)}$, i.e. these two brackets span a Poisson pencil
\be
a\{\,,\,\}_{gl(N)}+b \{\,,\,\}_r,\quad a,b \in \Bbb{C}.
\label{penc}
\ee

Moreover, if $r$ is standard, the bracket  $\{\,,\,\}_r$ (also called standard)
can be restricted to any  $GL(N)$-orbit in $gl(N)^*$. This property was proved in  \cite{D}.

Thus, the whole Poisson pencil (\ref{penc}) can be restricted to any $GL(N)$-orbit. As the Ki\-ril\-lov-Ko\-stant-Sou\-riau bracket is symplectic
on such an orbit, it is possible to use the method of the paper \cite{F} in order to construct integrable systems. We plan to go back to the
problem of constructing such systems in our subsequent publications. In the current paper we only show (by elementary methods) that the
bracket $\{\,,\,\}_r$ can be restricted to any generic orbit in $gl(N)^*$. Also, we consider the corresponding quantum algebra, called 
{\em the braided generic orbit}.

The paper is organized as follows. In the next section we introduce our basic objects: Hecke symmetries and RE algebras.  In section 3 we exhibit certain central elements of the RE
algebras and describe  relations between these elements. In section 4 we compare two approaches to constructing $\LL(R)$-modules similar to irreducible $U(gl(N))$-ones. In section
5 we introduce the GZ algebras, related to the RE algebras and their "restrictions" to {\em braided generic orbits}. In the last section we consider the semiclassical  structures corresponding to the above braided algebras.

 \section{Symmetries and RE algebras: basic definitions and examples}

By  a braiding we mean an operator $R:\vv\to \vv$ meeting the braid relation
$$
(R\ot I)(I\ot R)(R\ot I)=(I\ot R)(R\ot I)(I\ot R).
$$
Hereafter, $V$ is a finite dimensional vector space over the ground field $\C$, $\dim_{\,\C} V=N$, and $I$ stands for the identity matrix or operator.

A braiding $R$ satisfying the quadratic relation
\be
(q I-R)(\qq I+R)=0,\qquad q\in \C^\times,\quad q\not=\pm 1,
\label{Hecke}\
\ee
is called a {\it Hecke symmetry}. A complex number $q$ is assumed to be {\it genereic}: $q^k\not=1$ $\forall \,k\in{\Bbb Z}_+$.
As a consequence, a $q$-number $k_q = (q^k-q^{-k})/(q-q^{-1})\not=0$ for any nonzero integer $k$. 

If a braiding $R$ meets the relation $R^2=I$, it is called an {\it involutive symmetry}. By fixing a basis $\{x_i\}$ in the space $V$,
we identify $R$ with its matrix in the basis $\{x_i\ot x_j\}$ of $V^{\otimes 2}$.

The Drinfeld-Jimbo Hecke symmetries are defined as follows
\be
R=q\sum_{i =1}^N E_{i}^{i}\ot E_{i}^{i}+\sum_{i\not= j}^N E_{i}^{j}\ot E_{j}^{i}  + (q-\qq)\,\sum_{i< j}\,E_{i}^{i}\ot E_{j}^{j},
\label{He}
\ee
where $E_i^j$ is a matrix with the only nontrivial entry, which equals to 1, at the position $(i,j)$. If $N=2$, the corresponding Hecke symmetry is the first of the following two
matrices (hereafter we put $\la=q-\qq$)
\be  \left(\begin{array}{cccc}
q&0&0&0\\
0&\la&1&0\\
0&1&0&0\\
0&0&0&q\\
\end{array}\right),\qquad
\left(\begin{array}{cccc}
q&0&0&0\\
0&\la&1&0\\
0&1&0&0\\
0&0&0&-\qq\\
\end{array}\right). \label{exa} \ee
The second matrix  is a deformation of the super-flip\footnote{The notation $P_{m|n}$ stands for the super-flip in the space $\vv$, where $V=V_0\oplus V_1$, and
$V_0$ (resp. $V_1$) is the even (resp., odd) component with dimension $m$ (resp., $n$).}  $P_{1|1}$.

Also, there are known the so-called Cremmer-Gervais Hecke symmetries \cite{CrG}. Let us exhibit an example of such a symmetry for $N=3$:
\be \left(\begin{array}{ccccccccc}
q  &  0  &  0  &  0  & 0 & 0 & 0 & 0 & 0\\
\rule{0pt}{4.5mm}
0 & \la & 0 & \al & 0 & 0 & 0 & 0 & 0\\
\rule{0pt}{4.5mm}
0 & 0 & \la & 0 & 0 & 0 & q \,\al^2 & 0 & 0\\
\rule{0pt}{4.5mm}
0 & \al^{-1} & 0 & 0 & 0 & 0 & 0 & 0 & 0\\
\rule{0pt}{4.5mm}
0 & 0 & -q^{-3}\beta \al^{-1} & 0 & q & 0 & \beta & 0 & 0\\
\rule{0pt}{4.5mm}
0 & 0 & 0 & 0 & 0 & \la & 0 & \al & 0\\
\rule{0pt}{4.5mm}
0 & 0 & q^{-1}\al^{-2} & 0 & 0 & 0 & 0 & 0 & 0\\
\rule{0pt}{4.5mm}
0 & 0 & 0 & 0 & 0 &  \al^{-1}  & 0 & 0 & 0\\
\rule{0pt}{4.5mm}
0 & 0 & 0 & 0 & 0 & 0 & 0 & 0 & q
\end{array}\right).
\label{CG} \ee
Here $\al\not=0$ and $\beta$ are arbitrary complex numbers.

New examples of Hecke symmetries can be constructed by means of the following claim.
\begin{proposition}
\label{prop:1}
Let $V_i$, $i=1,2$, be two finite dimensional vector spaces and $R(i): \vv_i\to \vv_i$ be two Hecke symmetries. Let us define an operator
$$
R:(V_1\oplus V_2)^{\ot 2} \to (V_1\oplus V_2)^{\ot 2}
$$
by setting $R=R(i)$ on $\vv_i$ and
$$
R(x\ot y)=\la x\ot y+\al\, y\ot x,\quad R(y\ot x)=\al^{-1} x\ot y \quad \forall \,  x\in V_1,\,\, \forall \, y\in V_2,
$$
where $\al$ is an arbitrary nonzero  complex number.  Then $R$ is a Hecke symmetry.
\end{proposition}

Note that the method described in the above Proposition (in a slightly different form) was suggested in \cite{G}. Later it was generalized in \cite{MM} and  called {\em glueing method}.
On setting $q=1$ in Proposition \ref{prop:1}, we get a way of glueing involutive symmetries.

On taking  $\dim V_1=\dim V_2=1$ and setting $R(1)=R(2)=q$, we get the first Hecke symmetry from (\ref{exa}). If we put $R(2)=-\qq$, we get the second Hecke symmetry from
(\ref{exa}). We can continue this glueing procedure\footnote{At each step we can take different  $\al\not=0$.} and get Hecke symmetries,
which are called {\em almost standard} or {\em almost super-standard} if $R$ are not even. The reader is referred to \cite{GPS3} for the
notions of an even symmetry, the bi-rank of a symmetry and others. Emphasize, that when constructing the Hecke 
symmetries by above glueing procedure, we do not need any Hopf algebras (see \cite{R} where such Hopf algebras were constructed).

Below, we are  dealing with skew-invertible Hecke symmetries. Recall that  a braiding $R$ is called {\em skew-invertible} if there exists an operator $\Psi:\vv\to \vv$ such that
\be
\Tr_{(2)} R_{12}\Psi_{23}=P_{13}\quad \Leftrightarrow \quad \sum_{a,b =1}^N R_{ib}^{ja} \Psi_{ak}^{b\,n}=\de_i^n\de_k^j.
\label{Psi}
\ee

Given a skew-invertible braiding $R$, we define the corresponding $R$-trace $\Trr $ by setting
\be
\Trr X=\Tr (C X), \quad  C_1=\Tr_{(2)} \Psi_{12} \quad \Leftrightarrow \quad C_i^j=\sum_{k=1}^N\Psi_{ik}^{jk}
\label{CC}
\ee
for any $N\times N$ matrix $X$ (may be with noncommutative entries). The symbol $\Tr$  stands for the usual trace.

A Reflection Equation (RE) algebra $\LL(R)$ (without a spectral parameter) is a unital associative algebra over the complex field $\C$
generated by entries of the matrix
$L=\|l_i^j\|_{1\leq i,j \leq N}$ subject to the following matrix equality
\be
R\, (L\ot I)\, R\, (L\ot I)-(L\ot I)\, R\, (L\ot I)\, R=0.
\label{RE}
\ee

Note that if  an involutive or Hecke symmetry $R$ is a deformation of the usual flip $P$, the corresponding algebra $\LL(R)$ is a deformation of the the commutative algebra
$\Sym(gl(N))$. If $R$ is a Hecke symmetry,  this property means that for a generic $q$ the dimensions of the homogenous components of the algebra $\LL(R)$ are
classical\footnote{Note that for a Lie algebra $g$, belonging to one of the series $B_N,C_N, D_N$, the  algebra $\LL(R)$ with $R$ coming
from the corresponding quantum group, is {\it not} a deformation of $U(g)$. Nevertheless, some quotient of the algebra $\LL(R)$ is a
deformation of the function algebra $Fun(G)$ on the corresponding group $G$.}.

Observe that a linear shift of generators by the unit element $e_{\!\mbox{\tiny $\cal L$}}$
\be
l_i^j= e_{\!\mbox{\tiny $\cal L$}}\,\de_i^j-(q-\qq) \hl_i^j\quad \Leftrightarrow \quad  L=I\,e_{\!\mbox{\tiny $\cal L$}}-(q-\qq) \hL
\label{cha}
\ee
leads to the quadratic-linear relations on the generators $\hl_i^j$:
\be
R\, (\hL\ot I)\, R\, (\hL\ot I) - (\hL\ot I)\, R\, (\hL\ot I)\,R = R\, (\hL\ot I) -  (\hL\ot I)\,R.
\label{mRE}
\ee

The  RE algebra expressed in terms of generators $\hl_i^j$ is called the {\em modified} RE algebra and is denoted $\hLL(R)$. If $R$ is a
Hecke symmetry, deforming the usual flip $P$,
the algebra $\hLL(R)$ tends to the enveloping algebra $U(gl(N))$ as $q\to 1$. If  $R$ is a deformation of a super-flip $P_{m|n}$, then the corresponding modified RE algebra
tends to the enveloping algebra $U(gl(m|n))$. In general, for any skew-invertible  Hecke or involutive symmetry $R$ the modified RE algebra $\hLL(R)$ can be treated as the enveloping algebra of a "braided Lie algebra" (see \cite{GPS3}). Observe that it is not a deformation in the class of the enveloping algebras of  usual Lie algebras even though $R$ is a deformation
of the usual flip. Also, note that if $R$ is an involutive Hecke symmetry, the corresponding algebras $\LL(R)$ and $\hLL(R)$ are not
isomorphic  to each other (for instance, such as the algebras $\LL(P) = \Sym(gl(N))$ and $\hLL(P)= U(gl(N))$).

Thus,  we have two forms of the RE algebra --- modified and not --- and consequently two forms of braided analogs of the chain (\ref{chain}),
provided the initial Hecke symmetry is suitably chosen (see Lemma \ref{lem:2} below).  Constructing these braided analogs is the
main objective of the present paper.

\begin{lemma}
\label{lem:2} \rm
Let  $R=\|R_{ij}^{kl}\|_{1\leq i,j \leq N}$ be a Hecke symmetry such that it admits a sub-matrix $R'=\|R_{ij}^{kl}\|_{1\leq i,j,k,l \leq M}$, $M< N$, which is also a Hecke symmetry. Let, in addition, the following property takes place: if {\it any pair} of indices $i,j,k,l$ is less or equal to $M$ and at least one of the other indices is greater than $M$, then $R_{ij}^{kl}=0$.
Then the RE algebra $\LL(R')$ coincides with the subalgebra in $\LL(R)$ generated by elements $l_i^j$, $1\leq i,j \leq M$. A similar statement is valid for the corresponding modified
RE algebras $\hLL(R)$ and $\hLL(R')$.
\end{lemma}

If in this lemma we can take $M=N-1$ and if this procedure can be continued,  we get a longest possible chain of the RE algebras (modified or not) embedded into one another.

Note that the almost (super-)standard Hecke symmetries, introduced above, meet the condition of this lemma whereas the Cremmer-Gervais
symmetry (\ref{CG}) does not.
Indeed even though the components $R_{ij}^{kl}$ with $i,j,k,l\leq 2$ form a Hecke  symmetry, the entry $R_{22}^{13}$ does not vanish provided $\beta\not=0$.

\section{Some central elements of RE algebras and their properties}

Now, we exhibit a regular way of constructing  some central elements of any  algebra $\LL(R)$.  Let us  introduce  the following notation:
$$
R_k = I^{\otimes (k-1)}\otimes R\otimes I^{\otimes (p-k-1)}\in (\mathrm{Mat}_N(\C))^{\otimes p}, \quad 1\le k\le p-1.
$$
This is an embedding of $R$ into the space of $N^p\times N^p$ matrices, $p\ge 2$. Below we do not specify the concrete value of an integer
$p$ just assuming it to be large enough so that all the matrix formulas make sence.

Introduce also analogous embeddings for the generating matrix $L$ of the RE algebra $\LL(R)$:
\be
L_{\ov 1}=L_{\underline{1}} = L\otimes I^{\otimes (p-1)},\quad L_{\ov{k+1}} = R_{k}\, L_{\ov k}\, \RR_{k},
\quad L_{\underline{k+1}} = R^{-1}_{k}\, L_{\underline k}\, R_{k},\quad 1\le k\le p-1.
\label{R-copy}
\ee
Note, that all the matrices $R_m$, $L_{\ov n}$ and $L_{\underline n}$ are of the same size $N^p\times N^p$ for any $m$ and $n$.
Let us also denote
\be
L_{\ov{1\to k}}=L_{\ov 1}\, L_{\ov 2}\dots L_{\ov k}.
\label{string}
\ee
The following claim holds true.

\begin{proposition} {\rm \cite{IP}} Let $f(R_{ 1},\dots , R_{{k-1}})$ be an arbitrary polynomial in matrices $R_i$, $1\le i\le k-1$.
Then the element
\be
p_f (L) = \Tr_{R(1\dots k)}\left(f(R_{1},\dots,R_{{k-1}})\, L_{\ov{1\to k}}\,\right) = \Tr_{R(1\dots k)} \left( L_{\ov{1\to k}}\,\,f(R_{ 1},\dots, R_{{k-1}})\right)
\label{sem}
\ee
is central in the algebra $\LL(R)$.
\end{proposition}
Hereafter, we use the notation: 
$$
\Tr_{R(1\dots k)}(X)=\Tr_{R(1)}(\dots (\Tr_{R(k)}(X))\dots).
$$
Observe, that the second equality in (\ref{sem}) is valid in virtue of the cyclic property of the $R$-trace, which takes place for any polynomial in  $R_i$, $1\le i\le k-1$. This cyclic
property for the $R$-trace is a direct consequence of the relations $R_i \,C_i\, C_{i+1}=C_i\, C_{i+1} R_i$.

It is not difficult to show that the elements (\ref{sem}) form an associative algebra. Following \cite{IP}, we call this algebra {\em characteristic}. Conjecturally, the characteristic algebra
coincides with the centre $Z(\LL(R))$ but it is not yet shown  that any element of this center can be presented in the form (\ref{sem}).

Let us exhibit some central elements of special interest.  First, these are the elements
\be
p_k(L)= \Tr_{R(1\dots k)} (L_{\ov{1\to k}}\,R_{{k-1}}R_{{k-2}}\dots R_{1})
\label{pows}
\ee
which are called  the {\em power sums}. Note that they can be  reduced to the form $p_k(L)=\Tr_R L^k$.

Another important family of central elements consists of the so-called Schur polynomials  $s_\la(L)$, which are labeled by the partitions  $\la=(\la_1\geq\dots \geq \la_N)$.
Their explicit construction is given in \cite{GPS2}. A remarkable property of the Schur polynomials $s_\lambda(L)$ is their pure classical multiplication rule:
$$
s_\la(L)s_\mu(L)=\sum_\nu C_{\la\, \mu}^\nu\, s_\nu(L),
$$
where the coefficients  $C_{\la\, \mu}^\nu$ coincide with the Littlewood-Richardson ones (see \cite{GPS2} for a detailed proof). Consequently,
the Jacobi-Trudi formulae are also valid. These facts justify the term ``Schur polynomials'' for $s_\lambda(L)$.

Two particular cases of the Schur polynomials, namely the complete $h_k(L)$ and elementary $e_k(L)$ symmetric polynomials
are obtained by means of formula  (\ref{sem}) where the polynomial $f$ should be replaced respectively by the $R$-symmetrizer or by the
$R$-skew-symmetrizer in the space $V^{\ot k}$ (see \cite{GPS2}). The complete and elementary symmetric
polynomials are connected by the classical Wronski identities:
$$
\sum_{r=0}^k(-1)^r\,h_r(L)\,e_{k-r}(L) = 0\qquad \forall\, k\ge 1.
$$

The following  property of the RE algebra $\LL(R)$  is very important for our study. The generating  matrix $L$ of this algebra is subject 
to the Cayley-Hamilton
identity with central coefficients. If $R$ is even and its bi-rank is $(m|0)$ (see \cite{GPS3}),  the Cayley-Hamilton identity reads:
\be
L^m-q \,e_{1}(L)\, L^{m-1}+q^2\, e_{2}(L)\, L^{m-2}+\dots+ (-q)^{m-1}\, e_{m-1}(L)\, L+ (-q)^{m}\, e_{m}(L)\, I=0.
\label{CHH}
\ee

So, it is natural to introduce the eigenvalues $\mu_i$, $1\le i\le m$, as solutions of the system of polynomial equations:
\be
\sum_{i=1}^m \, \mu_i=q \,e_{1}(L),\qquad \sum_{1\le i_1<i_2\le m}\mu_{i_1}\mu_{i_2} = q^2e_2(L), \qquad \dots\qquad  \prod_{i=1}^m\, \mu_i=q^m\, e_{m}(L).
\label{q-e}
\ee
Thus, the eigenvalues $\mu_i$ can be treated as elements of the algebraic extension of the center of the RE algebra. Below we deal with the extended algebra
$\LL(R)[\mu_1,\dots,\mu_m]$.

The power sums are expressed in terms of spectral values by the following formula:
\be
p_k(L)=q^{-1}\,  \sum_{i=1}^m\, \mu_i^k\, \prod_{j\not= i}^m \, \frac{\mu_i-q^{-2}\, \mu_j}{ \mu_i-\mu_j}.
\label{pmu}
\ee

If a Hecke symmetry $R$ is of the general bi-rank $(m|n)$, then the corresponding Cay\-ley-Ha\-mil\-ton polynomial is of degree $m+n$ and it is not unital.
In this case the eigenvalues are naturally  split in two subfamilies: ``even'' eigenvalues $\mu_i$, $1\leq i \leq m$, and ``odd'' ones $\nu_j$, $1\leq j \leq n$.
They are also assumed to be central in the extended algebra $\LL(R)[\mu_1\dots \mu_m, \nu_1\dots\nu_n]$
 (see \cite{GPS2, GPS4} for detail).

Let us emphasize that for a skew-invertible Hecke symmetry of any bi-rank the power sums, the complete and elementary 
symmetric polynomials can be expressed via each other as follows
 $$
p_k=\det \left(\begin{array}{ccccc}
1_q\,e_1&1&0&...&0\\
2_q\,e_2&q\, e_1&1&...&0\\
3_q\,e_3&q^2\, e_2&q\, e_1&...&0\\
...&...&...&...&1\\
k_q\, e_k&q^{k-1}\, e_{k-1}&q^{k-2}\, e_{k-2}&...&q\, e_1
\end{array}\right),
$$
$$
p_k=(-1)^{k-1}\,\det \left(\begin{array}{ccccc}
1_q\,h_1&1&0&...&0\\
2_q\,h_2&\qq\, h_1&1&...&0\\
3_q\,h_3&q^{-2}\, h_2&\qq\, h_1&...&0\\
...&...&...&...&1\\
k_q\, h_k&q^{-(k-1)}\, h_{k-1}&q^{-(k-2)}\, h_{k-2}&...&\qq\, h_1
\end{array}\right),
$$
$$
 k_q !\, e_k=\det
\left(\begin{array}{ccccc}
p_1&1_q&0&...&0\\
p_2&q\, p_1&2_q&...&0\\
p_3&q\, p_2&q^2\, p_1&...&0\\
...&...&...&...&(k-1)_q\\
p_k&q\,p_{k-1}&q^2\, p_{k-2}&...&q^{k-1}\, p_{1}
\end{array}\right).
$$
$$
k_q !\, h_k=\det
\left(\begin{array}{ccccc}
p_1&-1_q&0&...&0\\
p_2&\qq\, p_1&-2_q&...&0\\
p_3&\qq\, p_2&q^{-2}, p_1&...&0\\
...&...&...&...&-(k-1)_q\\
p_k&\qq\,p_{k-1}&q^{-2}\, p_{k-2}&...&q^{-(k-1)}\, p_{1}
\end{array}\right).
$$

\begin{remark} \rm The q-analogs of power sums as functions in commutative variables have been introduced
 in the paper \cite{WY}, where  similar determinant formulae have been also exhibited. Whereas our determinant
relations can be treated in two ways, since all their elements can be expressed in terms of the
generators $l_i^j$ and in terms the eigenvalues $\mu_i$ as well.
\end{remark}

\section{Representations of RE algebras: two approaches}

There exist two approaches to constructing finite dimensional modules over the RE algebras. One of them is applicable to any RE algebra
$\LL(R)$  corresponding to an arbitrary skew-invertible Hecke symmetry $R$. The other one is working only if $R$ is standard.

The latter approach  is based on the seminal paper \cite{FRT}.  As follows from construction of this paper, the RE algebra $\LL(R)$, corresponding to the standard
Hecke symmetry $R$, can be embedded into $U_q(gl(N))$ as a subalgebra. Thus, the algebra $\LL(R)$ (and consequently $\hLL(R)$) becomes automatically
represented in any $U_q(gl(N))$-module. The category of all finite dimensional irreducible $U_q(gl(N))$-modules for a generic $q$ is well-known and is a deformation
of the category  of $U(gl(N))$-modules. So, with any pattern $\mm$ one can associate a $U_q(gl(N))$-module $\vm^q$, which coincides with the $U(gl(N))$-module 
$\vm$ as a vector space. Moreover, the action of the modified RE algebra $\hLL(R)$ in $U_q(gl(N))$-module $\vm^q$ is a deformation of the $U(gl(N))$ action on $\vm$.
Nevertheless, for a general $R$  this method is not applicable since the corresponding quantum group has no explicit description. 

Now, we briefly describe another way of $\LL(R)$-modules construction, valid for any algebra $\LL(R)$ whatever a skew-invertible Hecke
symmetry $R$ is (for a more detailed exposition of the construction the reader is referred to \cite{GPS3}).

This construction is performed in two steps. At the first step we construct a linear action of the algebra $\LL(R)$ on the spaces $V^{\ot n}$ for any $n\in {\Bbb Z}_+$. On the elements
of the tensor basis of $V^{\otimes n}$ this action is
defined by the rule:
\be
 L_{\un{n+1}}\triangleright(x_{|1\rangle }\ot \dots \ot x_{|n\rangle }) = J^{-1}_{n+1}\,x_{|1\rangle }\ot \dots \ot x_{|n\rangle },
\label{act}
\ee
where the symbol $\triangleright$ denotes the action of a linear operator and
$$
J_{{n+1}}= R_n\dots R_2\, R_1^2\, R_2\dots R_n
$$
is the image of the Jucys-Murphy element of the Hecke algebra $H_{n+1}(q)$ under its $R$-matrix representation (see \cite{GPS3} for detail).
At last, $x_{|r\rangle}\in V$ denotes an arbitrary vector of the basis set $\{x_i\}_{1\le i\le N}$ located at the position $1\le r\le n$ in the product $V^{\ot n}$.

At the second step we decompose the $\LL(R)$-module $V^{\ot n}$ into the direct sum of invariant submodules parameterized by all partitions $\mm\vdash n$
$$
V^{\otimes n} = \bigoplus_{\mm\,\vdash n}\bigoplus_{i=1}^{d_\mm} V_{\mm}^{(i)}
$$
where $d_{\mm}$ is the number of standard Young tables of the same shape $\mm$. The $\LL(R)$-submodules $V_{\mm}^{(i)}$ with the same $\mm$ and different $i$ are isomorphic
to each other. These submodules are images of the projection operators $P_{\mm}^{(i)}(R)$ corresponging to the primitive idempotents $e_\mm^{(i)}$ of the Hecke algebra $H_n(q)$
(see \cite{OP}) under the $R$-matrix representation:
\be
V_{\mm}^{(i)} = P_{\mm}^{(i)}(R)\triangleright V^{\otimes n}.
\label{Vm}
\ee

Now, we consider the connection of the representation (\ref{act}) with the representation of  the algebra $U(gl(N))$ in $V^{\otimes n}$. First, we pass to the generating matrix
$\hL$ by means of formula (\ref{cha}). Then, as directly follows from (\ref{act}), the matrix of the action of $\hL$ in the tensor basis $x_{|1\rangle}\dots x_{|n\rangle}$ reads:
\begin{equation}
\hL_{\underline{n+1}}  = \frac{I-J_{n+1}^{-1}}{q-q^{-1}} = R_{n}^{-1}+R_{n}^{-1}R_{n-1}^{-1}R_{n}^{-1}+\dots +R_{n}^{-1}\dots R_2^{-1}R_{1}^{-1}R_2^{-1}\dots R_{n}^{-1}.
\label{mrea-rep}
\end{equation}

Assuming the given Hecke symmetry $R$ to be a deformation of the usual flip $\lim_{q\rightarrow 1}R= P$ and by passing to the limit $q\rightarrow 1$ in (\ref{mRE}) we see that
the generators $\hl_i^j$ turn into those of $U(gl(N))$, while the matrix of representation (\ref{mrea-rep})  takes the form:
$$
\hL_{n+1} = P_{n}+P_{n}P_{n-1}P_{n}+\dots +P_{n}\dots P_{k+1}P_{k}P_{k+1}\dots P_{n}+\dots +P_{n}\dots P_2P_{1}P_2\dots P_{n}.
$$
Here, we used the fact that $\hL_{\un{n+1}}=\hL_{n+1}$ for $R=P$.

Taking into account that $P_{n}\dots P_{k+1}P_{k}P_{k+1}\dots P_{n} = P_{k \,n+1}$, where the notation $P_{ij}$ stands for the usual flip transposing the elements located at
the $i$-th and $j$-th positions,  we find that the $U(gl(N))$ generators $\hl_i^j$ are represented
in the tensor product $V^{\otimes n}$  by  the following matrix:
\be
\hL_{n+1} = P_{1\,n+1}+P_{2\, n+1}+\dots +P_{n\,n+1}.
\label{copr}
\ee
This formula means that the entries $\hl_i^j$ of the matrix $\hL_{n+1}$ are represented by the following operators in the space $V^{\otimes n}$:
$$
\hl_i^j \mapsto \sum_{m=1}^n I^{\otimes {m-1}}\otimes \hat E_i^j\otimes I^{\otimes {n-m}},
$$
where $\hat E_i^j$ is the linear operator corresponding to the matrix unit:
\be
\hat E_i^j\triangleright x_k = \delta_k^j \,x_i\quad \Leftrightarrow\quad \hL_2\triangleright x_{|1\rangle} = P_{12}x_{|1\rangle}.
\label{accc}
\ee

It should be emphasized that formula (\ref{copr}) is equivalent to the fact that the $U(gl(N))$-action (\ref{accc}) in the space $V$ is extended on the tensor product $V^{\otimes n}$
by multiple application of the usual coproduct defined in $U(gl(N))$.

As for the idempotents $P_\mm^{(i)}$, in the limit $q\to 1$ they turn into the corresponding idempotents of the group algebra $\C[{\cal S}_n]$ of the symmetric group ${\cal S}_n$
represented in $V^{\otimes n}$ by the usual flips. All these observations remain valid mutatis mutandis if the usual flip is replaced by another skew-invertible involutive symmetry.

Observe that this method of constructing the modules $\vm$ is implicitly based on the Schur-Weyl duality of the RE algebras and Hecke algebras described in \cite{GS3}.

Note that though the irreducibility of  the above $\LL(R)$-modules $V^{(i)}_\mm$ (\ref{Vm}) is not proved for a general Hecke symmetry $R$, the images of the central elements of the
algebra $\LL(R)[\mu_1\dots \mu_m]$ in these modules are scalar operators. Namely, as is shown in \cite{GSZ}, the images of the elements $\mu_k$ (up to the ordering of $\mu_k$)
in the module $V_\mm^{(i)}$, $\mm=(\mm_1\dots\mm_m)$ are scalar operators of the form:
\be
\mu_k\mapsto \chi_\mm(\mu_k)\,\mathrm{Id}_{V_\mm}, \qquad \chi_\mm(\mu_k) = q^{-2(\mm_k+m-k)}.
\label{chaa}
\ee
We emphasize, that the characters $\chi_\mm$ depend on the partition $\mm$ only and are the same for isomorphic modules $V_\mm^{(i)}$ with different values of $i$.

The formula (\ref{pmu}) allows us to compute the characters $\chi_\mm(p_k(L) )$ of the power sums. Namely, we have
\be
\chi_{\mm}(p_k(L))=q^{-2}\sum_{i=1}^m\, q^{-2k(\mm_i+m-i) }\,\prod_{j\not= i}^m \frac{(\mm_j-\mm_i+i-j+1)_q}{(\mm_j-\mm_i+i-j)_q}.
\label{mol}
\ee
Note that in the standard case formula (\ref{mol}) was obtained in \cite{JLM}. Also observe that in this case $m=N=\dim\, V$.

Using this formula we can compute the characters of the power sums $p_k(\hL)=\Trr \hL^k$ for the generating matrices $\hL$ of the algebras $\hLL(R)$.
Indeed, each matrix $\hL$ is also subject to a CH identity with central coefficients. Assuming $R$ to be even, we introduce the eigenvalues $\hmu_i$ of
the matrix $\hL$ in the same way as above. Formula (\ref{cha}) entails that the eigenvalues $\mu_i$ and $\hmu_i$ of the matrices $L$ and $\hL$ respectively
are related as follows
$$
\mu_i=1-(q-\qq)\hmu_i.
$$
This formula entails that the character of $\hmu_i$ on the module $V_\mm$ takes the value:
\be
\chi_\mm(\hmu_i)=\frac{1-q^{-2(\mm_i+m-i)}}{q-\qq}=q^{-(\mm_i+m-i)}(\mm_i+m-i)_q.
\label{chi}
\ee

From the other hand, formula (\ref{pmu}) allows one to express the power sums $p_k(\hL^k)=\Tr_R\hL^k$ in terms of the eigenvalues $\hmu_i$:
$$
p_k(\hL)=q^{-1}\,  \sum_{i=1}^m\, \hmu_i^k\, \prod_{j\not= i}^m \, \frac{\qq+\hmu_i-q^{-2}\, \hmu_j}{ \hmu_i-\hmu_j}.
$$
Assuming  $R$ to be a deformation of the usual flip and using formula (\ref{chi}), we get the formula from the paper \cite{PP} at the limit $q\to 1$.

As we noticed above, if $R$ is not even, the eigenvalues of the generating matrix $L$ split in two families $\{\mu_i\}$ and $\{\nu_i\}$. Though in this case a formula generalizing (\ref{pmu})
is known (see \cite{GPS4}), the values of $\mu_i$ and $\nu_j$ in the modules analogous to $\vm$ are not yet computed.

\section{GZ type algebras and braided orbits}

Now, we assume that a given Hecke symmetry $R$ meets the condition formulated in Lemma \ref{lem:2}. To be more concrete, we suppose
$R$ to be constructed by a series of glueings described after Proposition \ref{prop:1}, where at the first step we put $R(1)=R(2)=q$ with an
arbitrary $\al=\al_1\not=0$ and so on till the step number $m-1$. Afterwards, at the steps with numbers from $m$ to $m+n-1$ we put
$R(2)=-q^{-1}$. Thus, we get a Hecke symmetry of bi-rank $(m|n)$. Since the constructed Hecke symmetry meets the condition of Lemma
\ref{lem:2}, we get a chain of the RE algebras
\be
\LL_{(m+n-1)}(R) \supset \LL_{(m+n-2)}(R) \supset\dots \supset \LL_{(2)}(R)\supset  \LL_{(1)}(R)
\label{chainn}
\ee
as well as its modified version
\be
\hLL_{(m+n-1)}(R) \supset \hLL_{(m+n-2)}(R) \supset\dots \supset \hLL_{(2)}(R)\supset  \hLL_{(1)}(R)
\label{chainnn}
\ee
where the notation $\LL_{(k)}(R)$ stands for the RE algebra defined by the Hecke symmetry obtained at the $k$-th step of the above
glueing procedure.

Similarly to the classical setting or its super-analog we introduce commutative subalgebra $\ZZ_q$ as the union of the centers of all
components of the  chain (\ref{chainnn}). We call the algebra $\ZZ_q$ the {\em braided GZ algebra}.

Now, let us assume that all Hecke symmetries defining the RE subalgebras in (\ref{chainn}) are almost standard. Then the integers entering the exponents of formula (\ref{chaa})
and corresponding to the first algebra $\LL_{(N-1)}(R)$, are $\mm_1+N-1$,  $\mm_2+N-2$, etc. For the second algebra $\LL_{(N-2)}(R)$ we have $\mm_1'+N-2$,  $\mm_2'+N-3$,
and so on (here we use the notation accepted in (\ref{pat})). Note that these integers satisfy the interlacing condition. Consequently, it is so for the quantities (\ref{chaa})
if we assume the parameter $q$ to be real and $q>1$.  We consider this property as an analog of the interlacing property from \cite{GuS} which is valid for the eigenvalues of Hermitian
matrices, entering the construction of  the semiclassical counterpart of GZ model.

\begin{remark} \rm In Introduction we have described a way of constructing the GZ basis $e_{\cal M}$ in the $U(gl(N))$-module $\vm$. A method of constructing a similar basis of a
$U_q(gl(N))$-module $\vm^q$ was proposed in \cite{UTS}. This method is essentially based on the so-called lowering operators. Thus, 
if $R$ is standard, for the braided GZ algebra one
can construct the corresponding GZ basis, which  consists of quantum analogs $e_{\MM}^q$  of the vectors  $e_{\MM}$ considered in Introduction. The corresponding representation
of the braided GZ algebra can be realized via the embedding of the RE algebras into quantum groups at all levels. The characters of the
power sums of the generating matrices of all components $\LL(R_i)$ can be  computed according to formula (\ref{mol}).

Unfortunately, for a general $R$ we do not know how to decompose the $\LL_{(i)}(R)$-module into a direct sum of $\LL_{(i-1)}(R)$-modules,
since we do not know what the braided analogs of the lowering operators are in this case.
\end{remark}

Now, we assume $R$ to be an even Hecke symmetry of bi-rank $(m|0)$ and consider the problem of constructing braided analogs
of generic orbits in $gl(N)^*$. Let us recall that in the classical setting a generic orbit can be defined by the following system
\be
\Tr L^k= \sum_i\, \mu_i^k,\qquad  1\le k\le N.
\label{orb}
\ee
Here the matrix $L$ is the generating matrix  of the algebra $\Sym(gl(N))$ and its entries $l_i^j$ are considered to be linear functions on
$gl(N)^*$. The orbit defined by the system (\ref{orb}) is {\it generic}, iff it is the orbit of the matrix with pairwise distinct eigenvalues $\mu_i$.

Let us  consider the following quotient
\be
\OO=\LL(R)/\langle  \Trr L-\al_1(\mu),\, \Trr L^2-\al_2(\mu),\,\dots,\Trr L^m-\al_m(\mu)\rangle,
\label{ide}
\ee
where $\mu=(\mu_1,\mu_2,\dots,\mu_m)$ and $\al_k(\mu)$ are defined by the right hand side of (\ref{pmu}). We call this quotient a braided orbit.

Now, we are going to discuss the following problem: for which values of $\mu_i$ the braided orbit $\OO$ can be considered as generic?  According to the Serre-Swan approach, a classical orbit is a
regular variety iff the $\OO$-module of differential one-forms on it is projective. This criterion is also used in the noncommutative setting. We apply this criterion to our braided
orbits. We define this right $\OO$-module on such an orbit as follows
\be
\Omega^1(\OO)=({\Lambda}_R\ot \OO)/\langle d\Trr L,\,\dots\,, d\Trr L^m\rangle,
\label{diff}
\ee
 where $d$ is the de Rham operator and ${\Lambda}_R$ is a vector space spanned by the elements $d l_i^j$.

In order to avoid the problem of transposing the elements of the RE algebra and the differentials $d l_i^j$, we observe that in the classical setting, i.e. as $R=P$, in order to compute
$d\,\Tr L^k$ it suffices to apply  the de Rham operator $d$ only to the first factor of monomials entering $\Tr L^k$
$$
\Tr L^k \,\mapsto\,\Tr(dL)L^{k-1}, 
$$
since $\Tr$ operation on a matrix product is invariant under the cyclic permutations of the factors. So, the result of such an application of the de Rham operator differs from $d\,\Tr L^k$ by a factor $k$:
$$
d\,\Tr L^k = \Tr(dL)L^{k-1}+\Tr L(dL)L^{k-2}+\dots +\Tr L^{k-1}(dL) = k\,\Tr(dL)L^{k-1}.
$$

The cyclic invariance of the classical elements $\Tr L^k$ have the following braided analog.

\begin{lemma}\label{lem:6} \rm
The power sums (\ref{pows})  are invariant with respect to $R$-cyclic permutations defined on the products
$L_{\ov{1\to k}}$ by the rule:
$$
L_{\ov{1\to k}}\,\mapsto\, (R_{m-1}\, R_{m-2}\dots R_{1})^{-1}\,L_{\ov{1\to k}} \, R_{m-1}\, R_{m-2}\dots R_{1}.
$$
\end{lemma}

\noindent
{\bf Proof.} The claim of the lemma is a direct consequence of the definition of power sums and the cyclic property of the $R$-trace. \hfill\rule{6.5pt}{6.5pt}

\smallskip

For this reason, in the denominator of the quotient (\ref{diff}) we apply the de Rham operator $d$ only to the first factors in all summands of the elements $\Trr L^k$.

\begin{proposition}
The quotient {\rm (\ref{diff})} considered as a right $\LL(R)$-module is projective if the eigenvalues $\mu_i$, defining the braided orbit, are subject to the condition
\be
\mu_i\not= q^2\, \mu_j\qquad  1\le i,j\le m,\quad i\not=j.
\label{mmq}
\ee
\end{proposition}

This proposition was proved in \cite{GS1} under some conjecture concerning a property of the power sums. Now, we put forward no
conjecture using instead Lemma \ref{lem:6}.
Similarly to the mentioned conjecture this lemma allows one to apply the de Rham operator in the way described above. The remaining
part  of the proof from \cite{GS1} is unchanged.

Now, assume $R$ to be almost standard.
Let us fix a family of complex number $\mu_i$, $1\leq i \leq m$, subject to the condition (\ref{mmq}). Consider the corresponding quotient $\OO=\LL(R)/J_\mu$, where $J_\mu$ is the
ideal in the denominator of  (\ref{ide}). Then each term of the chain (\ref{chainn}) can be quotiented  over  the ideal $J_\mu$ with preserving all inclusions. Moreover, the central elements
of each component of the chain (\ref{chainn}) remain to be central in the corresponding quotient. The union of all centres of these quotients is also called a  {\em braided GZ algebra}.

If $R$ is almost super-standard and its rank is $(m|n)$, this construction remains valid mutatis mutandis. In this case the expressions of the  quantities $\al_i$, entering  formula
(\ref{ide}), become dependent on two families $\mu_i$ and $\nu_i$. Then the condition (\ref{mmq}) should be completed with a condition on $\nu_i$. The corresponding condition
is exhibited in \cite{GS1}.

The above construction of the braided GZ algebra on a braided generic orbit can be defined even if $R$ is almost super-standard. 
Also, it is possible to define braided generic orbits in  the corresponding modified RE algebras (see \cite{GS1}) and introduce 
related braided GZ algebras in the same way.

\section{Semiclassical counterpart of braided GZ algebras}

Let $R$ be a Hecke symmetry, deforming the usual flip $P$. Then by developing the matrix $P_{12} R_{12}$ in a series in $h$, where $h =
\log q$, we get
$$
P_{12} R_{12}=I+h\, r_{12}+o(h).
$$
The matrix $r$ is called the {\it classical $r$-matrix}.

It is well known and can be easily checked that $r$ satisfies the classical Yang-Baxter equation:
\be
 [r_{12}, r_{13}]+[r_{12}, r_{23}]+[r_{13}, r_{23}]=0.
\label{odin}
\ee
Also, the second degree equation (\ref{Hecke}) for the matrix $R$ entails that $r$ obeys the following relation:
\be
r_{12}+r_{21}=2\, P_{12},
\label{dva}
\ee
where $r_{21}=P_{12}\,r_{12}P_{12}$.

The linear in $h$ terms in the development of the RE algebra defining relation (\ref{RE}) leads to the following Poisson bracket on the commutative algebra $\Sym(gl(N))$
\be
\{L_1,\, L_2\}_r=r_{21}\, L_1\, L_2-L_1\, L_2\, r_{12}+ L_2\, r_{12}\, L_1- L_1\, r_{21}\, L_2.
\label{Poisson}
\ee
The fact that this expression  is a Poisson bracket follows from the relations (\ref{odin}) and (\ref{dva}). Indeed, due to the latter relation we get
$$
[r_{12}+r_{21}, L_1\, L_2]=0,
$$
since $L_1$ and $L_2$ commute with each other. This property implies the skew-symmetry of the bracket (\ref{Poisson}). The Jacobi identity for this bracket is ensured
by the classical Yang-Baxter equation (\ref{odin}).

Observe that for any classical $r$-matrix the corresponding bracket  $\{\,,\,\}_r$ has a very important property. It is compatible with the linear bracket $\{\,,\,\}_{gl(N)}$,
corresponding to the Lie $gl(N)$ structure. This fact can be checked by  straightforward computations. However, it also follows from the consideration of the corresponding
quantum (braided) object since the mentioned Poisson pencil is the semiclassical counterpart of the modified RE algebra (with a second parameter, introduced in front of the
linear term). Thus, we have a Poisson pencil (\ref{penc}).

Moreover, if $r$ is standard,  the bracket $\{\,,\,\}_r$  can be restricted to any semisimple orbit in $gl^*(N)$. This fact follows from the results 
of \cite{D}. (Note that the author
of \cite{D} dealt with the restriction of  the bracket $\{\,,\,\}_r$ to the space $sl(N)$ but the passage to the algebra $gl(N)$ is evident.)

However, if a given orbit is generic, this claim can be proved by the following elementary method.  First, recall that  any generic  orbit in
$gl^*(N)$ is obtained via  the system
(\ref{orb}) with pairwise distinct $\mu_i$.

\begin{proposition} The elements $\Tr L^k$ are Poisson central with respect to the bracket {\rm (\ref{Poisson})}
$$
\{l_i^j, \Tr L^k\}_r =0,\qquad \forall \, k\geq 0,\quad 1\leq i,j \leq N.
$$
\end{proposition}

\noindent
{\bf Proof.} For $k=1$ we have
$$
\{ \Tr L, L_2\}_r=\langle r_{21}\, L_1\rangle_1\, L_2-L_2 \,\langle L_1\, r_{12}\rangle_1+L_2 \,\langle  r_{12} \, L_1 \rangle_1 -\langle  L_1\, r_{21}\rangle_1\, L_2=0,
$$
where the notation $\langle\dots \rangle_1$ stands for applying the trace at the first position, and the cyclic property of the trace is used.

By using the Leibniz rule and the induction in $k$ it is not difficult to get the relation:
$$
\{L_1^k,\, L_2\}_r=r_{21}\, L_1^k\, L_2 - L_1^k\, r_{21}\, L_2 + L_2\, r_{12}\, L_1^k - L_2\, L_1^k\, r_{12} .
$$
Again, by applying the trace at the first position and using the cyclic property of the trace, we complete the proof.\hfill \rule{6.5pt}{6.5pt}

\smallskip

Thus, we  conclude that the  Poisson pencil (\ref{penc}) can be restricted to any generic orbit in $gl^*(N)$.

\begin{remark} \rm
Observe that the bracket $\{\,,\,\}_r$ can be restricted to the variety, defined by the equation $\det L=1$, where it coincides with  the
Semenov-Tian-Shansky bracket. (Recall that the Semenov-Tian-Shansky bracket is defined on the dual group (usually denoted $G^*$), 
where the initial group $G$ is
endowed with the Sklyanin bracket.) However, only on the whole algebra $\Sym(gl(N))$ it makes  sense to speak about its compatibility with  the linear Poisson bracket.

Also, note that there exists another way of constructing $r$-matrix brackets on the semisimple orbits, when $r$ is standard. It consists in
reducing the Sklyanin bracket from the group $SL(N)$ to a given orbit. As follows from the results of \cite{KRR}, this reduced bracket is
compatible with the Kirillov-Kostant-Souriau bracket iff the orbit is symmetric. The fact that the authors deal with  orbits of Hermitian 
matrices is not of matter.

Thus, on a symmetric orbit we have two Poisson pencils. It is possible to show that they coincide with each other.  Nevertheless, if a given
semisimple orbit is not symmetric only our  method  of constructing a Poisson pencil is valid.
\end{remark}

In order to construct GZ algebra in the Poisson algebra defined by the bracket $\{\,,\,\}_r$, it suffices to consider the generating matrices $L_{(k)}$ of the algebras $\Sym(gl(k))$,
defined as above, and to observe that the elements $\Tr L_{(k)}^p$ commute with each other with respect to this bracket and consequently to each bracket of the pencil (\ref{penc}).
Also, the same property is valid for the restrictions of this pencil to all semisimple orbits.

In \cite{F} this model was constructed for a symmetric orbit of $\Bbb{CP}^N$ type and some integrable systems on these orbits were
considered. We plan to consider similar integrable models on nonsymmetric (in particular, generic) orbits.

We complete the paper with the following example.  Consider the Hecke symmetry and the corresponding classical $r$-matrix of $GL(2)$ type:
$$
R = \left(\!\!
\begin{array}{cccc}
q&0&0&0\\
\rule{0pt}{4mm}
0&\lambda& x& 0\\
\rule{0pt}{4mm}
0 &x^{-1}&0&0\\
\rule{0pt}{4mm}
0&0&0&q\
\end{array}
\!\!\right),\qquad \,\,
r= \left(\!\!
\begin{array}{cccc}
1&0&0&0\\
\rule{0pt}{4mm}
0&-\alpha&0&0\\
\rule{0pt}{4mm}
0&2&\alpha&0\\
\rule{0pt}{4mm}
0&0&0&1\\
\end{array}\!\!\right),
$$
where $x$ is a function of $q$ with the following asymptotics at $q\rightarrow 1$: $x = 1+h\alpha +o(|h|)$. Thus, this Hecke symmetry is almost standard. However, as can be
easily checked, the RE algebra $\LL(R)$ does not depend\footnote{It is not so for the so-called RTT algebra, which is the algebra of quantized functions
on the group $GL(2)$.} on $x$.

Now, let us  exhibit the bracket $\{\,,\,\}_r$ on the generators $a,b,c, d\in \Sym(gl(2))$
$$
\begin{array}{lcl}
\{a,b\}_r = -2 ab&\quad & \{b,c\}_r = 2a(d-a)\\
\rule{0pt}{5mm}
\{a,c\}_r = 2 ac&& \{b,d\}_r = - 2 ab\\
 \rule{0pt}{5mm}
\{a,d\}_r = 0 & &\{c,d\}_r = 2 ac.
\end{array}
$$

This bracket is compatible with the bracket $\{\,,\,\}_{gl(2)}$. The GZ algebra is generated by 3 elements $a,\, a+d,\, a^2+d^2+2 bc$.

We complete the paper with the following observation. As we said above, in the standard case the RE algebra can be embedded into 
the quantum group $U_q(gl(N))$. However, this quantum group does not give rise to any nontrivial Poisson structure since it is isomorphic
to the classical object, namely, to the enveloping algebra $U(gl(N))$. Only the coalgebraic structure of $U_q(gl(N))$ is deformed with 
respect to the classical object.  Whereas the standard RE algebra is a meaningful deformation of the algebra $\Sym(gl(N))$ and it gives 
rise to a semiclassical structure.


\begin{thebibliography}{GPS3}
	
\bibitem[CG]{CrG} Cremmer E., Gervais J.L., The Quantum Group Structure Associated
with Non-Linearly Extended Virasoro Algebras, Comm. Math. Phys., 134 (1990) 619--632.

\bibitem[D]{D} Donin J., Double quantization on coadjoint representations of simple Lie groups and its orbits,  arXiv:math/9909160

\bibitem[FRT]{FRT} Faddeev L., Reshetikhin N., Takhtadzhyan L., Quantization of Lie groups and Lie algebras, Leningrad Math. J. 1 (1990), 193--225.

\bibitem[F]{F} Foth P., Bruhat Poisson structure on $CP^n$ and integrable systems, J. Math. Phys. 43 (2002) 3124--3132.

\bibitem[GZ]{GZ}  Gelfand I., Zetlin M., Finite-dimensional representations of the group of unimodular
matrices, Dokl. Akad. Nauk SSSR 71 (1950), 825–828 (in Russian).

\bibitem[GuS]{GuS} Guillemin V.,  Sternberg S., The Gelfand–Cetlin system and quantization of the complex
flag manifolds, J. Funct. Anal. 52 (1983), 106--128.

\bibitem[G]{G} Gurevich D., Algebraic aspects of the Yang-Baxter equation,  Leningrad Math. J.  2 (1991) 801--828.

\bibitem[GPS1]{GPS1} Gurevich D., Pytov P., Saponov P., Braided Differential operators on Quantum algebras,  J. of Geometry and
Physics 61 (2011), pp. 1485--1501.

\bibitem[GPS2]{GPS2} Gurevich D., Pyatov P., Saponov P., Cayley-Hamilton theorem for quantum matrix algebras of $GL(m|n)$ type,
St Petersburg Math Journal  17, N 1 (2005), 157--179

\bibitem[GPS3]{GPS3} Gurevich D., Pyatov P., Saponov P., Representation theory of (modified) Reflection
Equation Algebra of the $GL(m|n)$ type,  Algebra and Analysis 20 (2008) 70--133 (English translation: St Petersburg Math. J.  20 (2009) 213--253.)

\bibitem[GPS4]{GPS4} Gurevich D., Pyatov P., Saponov P., Spectral parameterization for power sums of quantum supermatrices, Theor. and Math. Physics 159 (2009) 587--597.

\bibitem[GS1]{GS1} Gurevich D.,  Saponov P., Generic super-orbits in $gl^*(m|n)$ and their braided counterparts, J. of Geometry and Physics 60 (2010) 1411--1423.

\bibitem[GS2]{GS2} Gurevich D.,  Saponov P., From Reflection Equation Algebra to Braided Yangians, Proceedings of the 1st International Conference on Mathematical Physics,
Grozny, Russia, 2016. Springer Proceedings in Math. and Statistics V.273 (2018).

\bibitem[GS3]{GS3} Gurevich D.,  Saponov P., Quantum Schur-Weyl duality and q-Frobenius formula related to Reflection Equation algebras, International Math. Research Notices, 2025, no.3 (2025) rnae288. \par

\bibitem[GSZ]{GSZ} Gurevich D.,  Saponov P., Zaitsev M., Representations of Spectra of Quantum Matrices, in progress.

\bibitem[IP]{IP} Isaev A., Pyatov P., Spectral Extension of the Quantum Group Cotangent Bundle, Comm. Math. Phys. 288 (2009) 1137--1179.

\bibitem[JLM]{JLM} Jing N., Lui M., Molev A., Eigenvalues of quantum Gelfand invariants, J. Math. Phys. 65, 061703 (2024).

\bibitem[KW]{KW} Kostant B.,  Wallach N., Gelfand-Zetlin modules from the perspective of classical mechanics,  arXiv:math/0408342.

\bibitem[KRR]{KRR} Khoroshkin S., Radul A., Rubtsov V., A Family of Poisson Structures on Hermitian Symmetric Spaces, CMP 152 (1993) 299--315.

\bibitem[MM]{MM} Majid Sh., Markl M.,  Glueing operation for R-matrices, quantum groups and link-
invariants of Hecke type, Math. Proc. Cambridge Soc. 119 (1996) 139–166.

\bibitem[M]{M}  Molev A., Gelfand-Tsetlin bases for classical Lie algebras, arxiv.org/pdf/math/0211289.

\bibitem[OP]{OP} Ogievetsky O., Pyatov P., Lecture on Hecke algebras, preprint CPT-000/P.4076 (2004).

\bibitem[PP]{PP} Perelomov A., Popov V., Casimir operators for semisimple Lie lgebras, Izv. AN SSSR Ser. mat. 32 (1968), 1368--1390.

\bibitem[R]{R} Reshetikhin N., Quasitriangular Hopf algebras and invariants of links,  Leningrad Math. J., 1 (1990), 491–513

\bibitem[UTS]{UTS} Ueno K., Takebayashi T., Shibukawa Y., Gelfand-Zetlin basis for $U_q(gl(N+1))$-modules, Lett. Math. Phys. 18 (1989) 215--221.

\bibitem[WY]{WY} Wysbourne B., Yang M., $q$-Deformation of symmetric functions and Hecke algebras $H_n(q)$ of type $A_{n-1}$, AIP Conference Proceedings 266, 64 (1992).

\end{thebibliography}
\end{document}